\documentclass[leqno]{article}
\usepackage{amsmath,amsthm,amsfonts,amscd,eucal}

\numberwithin{equation}{section}

\def\cc{{\mathcal C}}

\def\cam{{\mathcal M}}
\def\cn{{\mathcal N}}

\def\cp{{\mathcal P}}

\def\bc{{\mathbb C}}

\def\bn{{\mathbb N}}
\def\br{{\mathbb R}}

\def\b{\beta}

\def\n{\nu}

\def\r{\rho}

\def\f{\varphi}
\def\c{\chi}
\def\itm#1{\item{$(#1)$}}

\def\gP{{\mathfrak P}}

\newtheorem{Thm}{Theorem}[section]
\newtheorem{Cor}{Corollary}[section]
\newtheorem{Prop}{Proposition}[section]
\newtheorem{Lem}{Lemma}[section]
\theoremstyle{definition}
\newtheorem{Def}{Definition}[section]
\newtheorem{Ex}{Example}[section]

\theoremstyle{remark}
\newtheorem{Rem}{Remark}[section] 
 
\begin{document}
\noindent

{\LARGE   
\textbf{\textsf{On the characterisation of paired monotone metrics}}} 

\bigskip

%{\LARGE   
%\textbf{\textsf{statistically monotone metrics}}} 

\bigskip

Paolo Gibilisco\footnote[1]{Electronic mail: gibilisc@sci.unich.it}

\textit{Dipartimento di Scienze, Facolt\`a di Economia, Universit\`a 
di 
Chieti-Pescara ``G. D'Annunzio'', Viale Pindaro 42, I--65127 Pescara, 
Italy.}

\bigskip

Tommaso Isola \footnote[2]{Electronic mail: isola@mat.uniroma2.it}

\textit{Dipartimento di Matematica, Universit\`a di Roma ``Tor 
Vergata'', Via della Ricerca Scientifica, I--00133 Roma, Italy.}

\markright{Paired monotone metrics}

 \begin{abstract}
     Hasegawa and Petz introduced the notion of paired monotone
     metrics.  They also gave a characterisation theorem showing that
     Wigner-Yanase-Dyson metrics are the only members of the paired
     family.  In this paper we show that the characterisation theorem
     holds true under hypotheses that are more general than those used
     in the above quoted references.
     
     Key words and phrases: Monotone metrics, Wigner-Yanase-Dyson information.
 \end{abstract}
  
 \section{Introduction} \label{sect.I}

 Monotone metrics are the quantum counterpart of Fisher
 information and are classified by Petz (1996,2002).  The
 Wigner-Yanase-Dyson information content (see Lieb (1973), Wigner and 
 Yanase (1963))
 $$
 I_p(\r,A)=- \hbox{Tr}([\r^p,A][\r^{1-p},A])
 $$
 can be seen as a one-parameter family of monotone metrics, see
 Hasegawa and Petz (1997).  There, Hasegawa and Petz gave a proof that the
 WYD-metrics are the only monotone metrics possessing a certain
 pairing property (in Hasegawa (2003) Hasegawa discusses how this reflects
 on the associated relative entropy along the lines of Lesniewski and 
 Ruskai (1999)). 
 This is substantially related to the pairing of the non-commutative versions of
 Amari embeddings
 $$
 \r \to \frac{\r^p}{p} \qquad \quad\r \to \hbox{log}(\r).
 $$
 The purpose of the present paper is to present a partially different
 proof of the characterisation theorem.  In Hasegawa and Petz (1997) 
 and Hasegawa (2003) a certain
 boundary behaviour is used as an hypothesis.

 Here we show that the characterisation theorem holds true under more
 general conditions, that is without the above hypothesis (see the
 Remark \ref{condizioni}).  While we use means that are relatively
 less elementary (the theory of regularly varying functions) it seems
 that the present proof also fills some gaps appearing in the
 arguments of Hasegawa and Petz (1997) 
 and Hasegawa (2003).  It should be emphasized that the
 pairing discussed here is related to the duality of non-commutative
 $\alpha$-connections as discussed in many papers (see Nagaoka (1995), 
 Hasegawa (1995), Gibilisco and
 Isola (1999), Amari and Nagaoka (2000), Grasselli and Streater
 (2001), Jen\v cova (2001), Grasselli (2002)).

 Another goal of this paper is to relate the above pairing to the
 duality of uniformly convex Banach spaces according to the lines of
 our previous works Gibilisco and Pistone (1998), Gibilisco and Isola
 (1999, 2001a,b): this appears, up to now, as one of the the main
 tools for the infinite dimensional approach to Information Geometry.
 
 The structure of the paper is as follows. In section 2 we present the 
 notion of pull-back of duality pairing and discuss the case of 
 commutative Amari embeddings. In section 3 we review the theory of 
  monotone metrics and their pairing. In section 4 one 
 finds the basic results on regularly varying functions that are needed 
 in the sequel. Section 5 contains the proof of the characterisation 
 theorem.

 \section{Pull-back of duality pairings} \label{sect.II}

 Let $V,W$ be vector spaces over $\br$ (or $\bc$).  One says that
 there is a duality pairing if there exists a separating bilinear form
 $$
 \langle \cdot, \cdot \rangle:V \times W \to \br.
 $$
 Let $\cam, \cn$ be differentiable manifolds.  A differentiable
 function $\f:\cam \to \cn$ is an immersion if its differential
 $D_{\r}\f:T_{\r}\cam \to T_{\f(\r)}\cn$ is injective, for any
 $\r\in\cam$.  

 \begin{Def}
     Suppose we have a pair of immersions $(\f,\c)$, 
     where $\f:\cam \to \cn$ and $\c: \cam
     \to {\tilde{\cn}}$,  such that a duality pairing exists between
     $T_{\f(\r)}\cn$ and $T_{\c(\r)}{\tilde{\cn}}$ for any $\r \in
     \cam$.  Then we may pull-back this pairing on $\cam$ defining
     $$
     \langle u,v \rangle_{\r}^{\f,\c}:=\langle D_{\r}\f(u),
     D_{\r}\c(u) \rangle \qquad u,v \in T_{\r}\cam.
     $$ 
 \end{Def}

 The most elementary example is given by the case where $\cn
 ={\tilde{\cn}}$ is a riemannian manifold, $\f =\c$ and the duality
 pairing is just given by the riemannian scalar product on
 $T_{\f(\r)}\cam$.  This is called the pull-back metric induced by the 
 map $\f$.

 A first non-trivial example is the following.  Let $X$ be a uniformly
 convex Banach space such that the dual $\tilde X$ is uniformly
 convex.  We denote by $\langle \cdot,\cdot \rangle$ the standard
 duality pairing between $X$ and $\tilde X$.  Let $J:X \to \tilde X$
 be the duality mapping, that is $J$ is the differential of the map $v
 \to \frac{1}{2}||v||^2$ (see Berger (1977), p.  373).  This implies
 that $J(v)$ is the unique element of the dual such that
 $$
 \langle v,J(v) \rangle= ||v||^2=||J(v)||^2.
 $$

 \begin{Def}
     Let $\cal M$ be a manifold.  If we have a map $\f: {\cal M} \to
     X$ we can consider a dualised pull-back that is a bilinear form
     defined on the tangent space of $\cal M$ by
     $$
     \langle A,B \rangle^{\f}_{\r}:=\langle A, B \rangle _{\r}^{\f,J
     \circ \f}=\langle D_{\r}\f(A),D_{\r}(J \circ \f)(B) \rangle.
     $$
 \end{Def}
 
 \begin{Rem}
     For $X$ a Hilbert space, $J$ is the identity, and this is again
     the definition of pull-back metric induced by the map $\f$.
 \end{Rem}
 
 \begin{Ex}\label{ex:1}
     Let $(X, {\cal F}, \mu)$ be a measure space.  If $f$ is a
     measurable function and $q \in (1,+\infty)$ then $||f||_q :=
     (\int|f|^q)^{\frac{1}{q}}$.  Moreover ${\tilde q}$ is defined by
     $\frac{1}{q}+\frac{1}{\tilde q}=1$.  Set
     $$
     L^q=L^q(X, {\cal F}, \mu)=\{ f \hbox{ is measurable}|\
     ||f||_q<\infty\}
     $$
     Define $N^q$ as $L^q$ with the norm
     $$
     ||f||_{N^q}:= \frac{||f||_q}{q}.
     $$
     Obviously $\widetilde{N^q}$ (the dual of $N^q$) can be
     identified with $N^{\tilde q}$.  Indeed if $f \in N^q$ and $g \in
     N^{\tilde q}$ define
     $$
     T_g(f):= \int \frac{f}{q} \frac{g}{\tilde q}
     $$
     One has
     $$
     ||T_g||= \hbox{sup} \frac{|T_g(f)|}{||f||_{N^q}}= \hbox{sup}
     \frac{\int \frac{f}{q} \frac{g}{\tilde q} }{\frac{||f||_q}{q}}=
     \frac{1}{\tilde q}\hbox{sup} \frac{\int fg}{||f||_q}=
     \frac{||g||_{\tilde q}}{\tilde q}= ||g||_{N^{\tilde q}}
     $$
     from this easily follows that $g \to T_g$ is an isometric
     isomorphism between $\widetilde{N^q}$ and $N^{\tilde q}$.  Now
     suppose that $\r>0$ is measurable and $ \int \r=1$, namely $\r$
     is a strictly positive density.  Then $v=q \r ^{\frac{1}{q}}$ is
     an element of the unit sphere of $N^q$ and it is easy to see that
     $J(v) = {\tilde q}\r^{\frac{1}{\tilde q}}$.  The family of
     maps $\r \to q\r^{\frac{1}{q}}$ are just the Amari embeddings.
 \end{Ex}
 
 Let $X=\{1,...,n\}$ and let $\mu$ be the counting measure.  In this
 case $N^q$ is just $\br^n$ with the norm $\frac{||\cdot||_q}{q}$. 
 Let ${\cal P}_n=\{v \in \br^n|v_i>0,\sum v_i=1\}$. 

 \begin{Prop}\label{Fisher}
     Consider the Amari embedding $\f:\r\in{\cal P}_n \to 
     q\r^{\frac{1}{q}}\in N^q$ for an
     arbitrary $q\in(1,+\infty)$.  Then the bilinear form
     $$
     \langle A,B \rangle^{\f}_{\r}:=\langle A, B \rangle _{\r}^{\f,J
     \circ \f}=\langle D_{\r}\f(A),D_{\r}(J \circ \f)(B) \rangle
     \qquad \qquad A,B \in T_{\r}{\cal P}_n
     $$
     is just the Fisher information.
 \end{Prop}
 \begin{proof}
     $$
     \langle D_{\r}\f(A),D_{\r}(J \circ \f)(B) \rangle = \int
     (\r^{\frac{1}{q}-1}A)(\r^{\frac{1}{\tilde q}-1}B)= \int
     \frac{AB}{\r}
     $$
\end{proof}

 The above result can be stated in much greater generality using the
 machinery of Pistone and Sempi (1995), Gibilisco and Isola (1999).
 
 \section{Paired monotone metrics} \label{sect.III}

 In the commutative case a Markov morphism  is a stochastic map
  $T: \br^n \to \br^k$.  In the noncommutative case a
 Markov morphism is a completely positive and trace preserving operator
 $T: M_n \to M_k$, where $M_n$ denotes the space of $n$ by $n$ complex
 matrices.  We shall denote by ${\cal D}_n$ the manifold of strictly
 positive elements of $M_n$ and by ${\cal D}^1_n \subset {\cal D}_n$
 the submanifold of density matrices.

 In the commutative case a monotone metric is a family of riemannian
 metrics $g=\{g^n\}$ on $\{\cp_n\}$, $n \in \bn$ such that
 $$
 g^m_{T(\rho)}(TX,TX) \leq g^n_{\rho}(X,X)
 $$
 holds for every stochastic map $T:\br^n \to \br^m$ and all $\rho
 \in \cp_n$ and $X \in T_\rho \cp_n$.

 In perfect analogy, a monotone metric in the noncommutative case is a
 family of riemannian metrics $g=\{g^n\}$ on $\{{\cal D}^1_n\}$, $n
 \in \bn$ such that
 $$
 g^m_{T(\rho)}(TX,TX) \leq g^n_{\rho}(X,X)
 $$
 holds for every stochastic map $T:M_n \to M_m$ and all $\rho \in
 {\cal D}^1_n$ and $X \in T_\rho {\cal D}^1_n$ (see Chentsov and Morotzova 
 (1990)).

 Let us recall that a function $f:(0,\infty)\to \br$ is called
 operator monotone if for any $n\in \bn$, any $A$, $B\in M_n$ such
 that $0\leq A\leq B$, the inequalities $0\leq f(A)\leq f(B)$ hold. 
 An operator monotone function is said symmetric if $f(x):=xf(x^{-1})$
 and normalised if $f(1)=1$.  With such operator monotone functions
 $f$ one associates the so-called Chentsov--Morotzova functions
 $$
 c_f(x,y):=\frac{1}{yf(xy^{-1})}\qquad\hbox{for}\qquad
 x,y>0.
 $$

 \begin{Prop} \label{cm}
     For a CM-function the following is true
     \itm{i} $c(tx,ty)=\frac{1}{t}c(x,y) \qquad \forall x,y,t>0$  
     \itm{ii} $c(x):= \lim_{y \to x} c(x,y)= \frac{1}{x}$.
 \end{Prop}

 Define $L_{\r}(A) := \r A$, and $R_{\r}(A) := A\r$.  Since $L_{\r}$
 and $R_{\r}$ commute we may define $c(L_{\r},R_{\r})$.  Now we can
 state the fundamental theorems about monotone metrics.  In what
 follows uniqueness and classification are stated up to scalars.

 \begin{Thm} (Chentsov 1982) 
     There exists a unique monotone metric on $\cp_n$ given by the
     Fisher information.
 \end{Thm}

 \begin{Thm} (Petz 1996) 
     There exists a bijective correspondence between monotone metrics
     on $M_n$ and normalised symmetric operator monotone
     functions.  This correspondence is given by the formula
     $$
     g_{\r}^{f} (A,B) :=Tr(A\cdot c_f(L_\rho,R_\rho)(B)).
     $$
 \end{Thm}

 The tangent space to ${\cal D}^1_{n}$ at $\r$ is given by
 $T_{\r}{\cal D}^1_{n} \equiv \{ A\in M_{n} : A=A^{*}, Tr(A)=0\}$, and
 can be decomposed as $T_{\r}{\cal D}^1_{n} = (T_{\r}{\cal
 D}^1_{n})^{c} \oplus (T_{\r}{\cal D}^1_{n})^{o}$, where $(T_{\r}{\cal
 D}^1_{n})^{c}:= \{ A\in T_{\r}{\cal D}^1_{n} : [A,\r] = 0\}$, and
 $(T_{\r}{\cal D}^1_{n})^{o}$ is the orthogonal complement of
 $(T_{\r}{\cal D}^1_{n})^{c}$, with respect to the Hilbert-Schmidt
 scalar product $\left< A,B \right>_{HS} := Tr(A^{*}B)$.  A typical element
 of $(T_{\r}D_{n})^{o}$ has the form $A=i[\r,U]$, where $U$ is
 self-adjoint.  Each statistically monotone metric has a unique
 expression (up to a constant) given by $Tr(\r^{-1}A^{2})$, for $A\in
 (T_{\r}{\cal D}^1_{n})^{c}$.
 
 \begin{Prop}\label{derivative} (See Bhatia 1997) 
     Let $A \in T_{\r}{\cal D}^1_{n}$ be decomposed as $A=A^c+i[\r,U]$
     where $A^c \in (T_{\r}{\cal D}^1_{n})^{c}$ and $i[\r,U] \in
     (T_{\r}{\cal D}^1_{n})^{o}$.  Suppose $\f \in {\cal
     C}^1(0,+\infty)$.  Then
     $$
     (D_{\r}\f)(A)=\f'(\r)A^c+i[\f(\r),U].
     $$
 \end{Prop}

      Let $\f,\c \in {\cal C}^1(0,+\infty)$.  Using the functional
      calculus one may consider $(\f,\c)$ as a pair of functions from
      ${\cal D}^1_{n}$ to $M_n$, for which a duality pairing is 
      provided by the Hilbert Schmidt scalar product $\left<\cdot , 
      \cdot \right>_{HS}$. Therefore, according to the previous section, we 
      can define the paired metric induced by $(\f,\c)$ as
      $$
      \langle A,B \rangle_{\r}^{\f,\c}=\hbox{Tr}(D_{\r}\f 
      (A) \cdot D_{\r}\c(B)), \qquad A,B \in T_{\r}{\cal D}_{n}^{1}. 
      $$

 \begin{Prop}\label{dualcm} (Hasegawa and Petz (1997), Hasegawa 
 (2003))
     Let $f$ be operator monotone, $c=c_{f}$ the associated
     CM-function.  For a pair $\f,\c \in {\cal C}^1(0,+\infty)$, the
     equality
     $$
      \langle A,B \rangle_{\r}^{\f,\c}=\hbox{Tr}(A\cdot 
      c(L_{\r},R_{\r})(B)).
     $$
     implies
     \begin{equation}\label{eq:dualpair}
	 c(x,y)=\frac{\f(x)-\f(y)}{x-y}\cdot \frac{\c(x)-\c(y)}{x-y}.
     \end{equation}
 \end{Prop}
 \begin{proof}
     It is enough to consider elements of $(T_{\r}D_{n})^{o}$. 
     Suppose $A=i[\r,U]$ and $B=i[\r,V]$ where $U,V$ are self-adjoint. 
     Using Proposition \ref{derivative} one has
     $D_{\r}\f(A)=i[\f(\r),U]$ and similarly for $B$.  Therefore
     $$
     \langle A,B \rangle_{\r}^{\f,\c}=\langle D_{\r}\f (A),D_{\r}\c(B)
     \rangle= \langle i[\f(\r),U], i [\c(\r),V] \rangle = \langle
     i\hat\f(L_{\r},R_{\r})[\r,U], i \hat\c(L_{\r},R_{\r})[\r,V]
     \rangle,
     $$
     where $\hat\f(x,y) := \frac{\f(x)-\f(y)}{x-y}$, and similarly 
     for $\hat\c$.  On the other
     hand it is true that
     $$
     \hbox{Tr}(A\cdot c(L_\rho,R_\rho)(B))=
     \hbox{Tr}(i[\r,U]c(L_{\r},R_{\r})(i[\r,V])).
     $$
     From the above equations and the arbitrariness of $A,B$ one has
     the conclusion.
 \end{proof}
 
 \begin{Def}
     In the hypotheses of Proposition \ref{dualcm}, we say that $\langle 
     \cdot, \cdot \rangle_{\r}^{\f,\c}$ is a {\it paired monotone metric}.
     Moreover we set
     $$
     \gP:=\{(\f,\c) | \f,\c \in {\cal C}^1(0,+\infty) \hbox{ and } \f,\c 
     \hbox{ induce a paired monotone metric} \}
     $$
 \end{Def}
 
 In what follows we give examples of elements of $\gP$.

 \begin{Def}
     $$
     f_p(x):=p(1-p) \frac{(x-1)^2}{(x^p-1)(x^{1-p}-1)} \qquad 
     p\in\br\setminus\{0,1\}
     $$
     $$
     f_0(x)=f_1(x):=\frac{x-1}{\hbox{log}(x)}.
     $$
 \end{Def}
 
 Obviously $f_p=f_{1-p}$ and
 $$
 f_0=\lim_{p \to 0} f_p =\lim_{p \to 1} f_p=f_1.
 $$

 Moreover we have that $f_{-1}$ is the function of the RLD-metric,
 $f_0=f_1$ is the function of the BKM-metric and $f_{\frac{1}{2}}$ is
 the function of the WY-metric.

 \begin{Def}
     $$
     (\f_p(x),\c_p(x))=(\frac{x^p}{p},\frac{x^{1-p}}{1-p}) \qquad
     p\in \br\setminus\{0,1\}
     $$
     $$
     (\f_0(x),\c_0(x))=(\f_1(x),\c_1(x))=(x, \log x).
     $$
 \end{Def}
 
 \begin{Thm}\label{family}  {\rm Hasegawa and Petz (1997), Hasegawa 
 (2003)}
     $(\f_p,\c_p)$ induce a paired monotone metric if and only if $p \in [-1,2]$.
 \end{Thm}
 \begin{proof}
     The proof consists in showing that the function $f_p$ is operator
     monotone iff $p \in [-1,2]$.  
          
     After this one has immediately that 
     $$
     c_p(x,y)=\frac{1}{yf_p(\frac{x}{y})}=\frac{\f_p(x)-\f_p(y)}{x-y}\cdot
     \frac{\c_p(x)-\c_p(y)}{x-y}
     $$
     and this ends the proof.
 \end{proof}

 Now let $p\in(0,1)$ and set $q=\frac{1}{p}$.  We use again the
 symbol $N_q$ to denote $M_n$ with the norm
 $$
 ||A||_{N^q}=q^{-1}(\hbox{Tr}(|A|^q))^{\frac{1}{q}}
 $$
 All the commutative construction of Example \ref{ex:1} goes through. 
 The following Proposition is the non-commutative analogue of
 Proposition \ref{Fisher} (see also Hasegawa and Petz (1997), Jen\v 
 cova (2001), Gibilisco and Isola (2001b), 
 Grasselli(2002)).

 \begin{Prop}\label{Banach} 
    Let $\f:\r\in {\cal D}^1_n \to q\r^{\frac{1}{q}}\in N_q$ be the
    Amari embedding.  The dualised pull-back
    $$
    \langle A,B \rangle^{\f}_{\r}:=\langle A, B \rangle _{\r}^{\f,J
    \circ \f}=\langle D_{\r}\f(A),D_{\r}(J \circ \f)(B) \rangle
    $$
    coincides with the Wigner-Yanase-Dyson information.
 \end{Prop}
 \begin{proof}
    It is a straightforward application of Proposition \ref{dualcm}.
 \end{proof}

\section{Regularly varying functions} \label{sect.IV}

For the content of this section see Bingham et al. (1987). 

 \begin{Def}
     Let $\ell$ be a measurable positive function defined on some
     neighbourhood $[X, +\infty)$ of infinity and satisfying
     $$
     \lim_{x \to + \infty} \frac{\ell(tx)}{\ell(x)}=1 \qquad \forall
     t>0;
     $$
     then $\ell$ is said slowly varying.
 \end{Def}
 
 \begin{Rem}
     Defining $\ell(x)=\ell(X)$ on $(0,X)$ one often considers $\ell$
     defined on $(0,+\infty)$.
 \end{Rem}
 
 Some examples of slowly varying functions are $\ell(x)=\log(x), \log
 (\log(x)), \hbox{exp}(\log (x)/ \log(\log(x)))$.

 \begin{Prop}\label{eq:limiti}
     If $\ell$ is slowly varying and $p>0$ then
     $$
     \lim_{x \to + \infty}x^p \ell (x)=+\infty \qquad \qquad \lim_{x
     \to + \infty} \frac{\ell(x)}{x^p}=0.
     $$ 
 \end{Prop}

 \begin{Def}
     A measurable function $h>0$ satisfying 
     $$
     \lim_{x \to + \infty} \frac{h(tx)}{h(x)}=t^p \qquad  \forall t>0
     $$
     is called regularly varying of index $p$; we write $h \in R_p$. 
     Therefore $R_0$ is the class of slowly varying functions.  We set
     $R:= \cup_{p \in \br}R_p$.
 \end{Def}
 
 \begin{Rem}
     Obviously homogeneous functions are very particular cases of
     regularly varying functions.
 \end{Rem}
 
 \begin{Prop}\label{RV}
     Assume $h>0$ is a measurable function, and there exists a
     function $j$ such that
     \begin{equation}\label{eq:j}
	 \lim_{x \to +\infty} \frac{h(tx)}{h(x)}= j(t) \in (0,+\infty)
     \end{equation}
     for all $t$ in a set of positive measure. Then
     \itm{i} the equation \ref{eq:j} holds for all $t>0$;
     \itm{ii} there exists $p\in \br$ such that $j(t)=t^p$, $\forall 
t>0$;
     \itm{iii} $h(x)=x^p \ell(x)$ where $\ell$ is slowly varying.
 \end{Prop}
 
 Sometimes, as in the present paper, one is interested in the
 behaviour at the origin.

 \begin{Def}
     If $h$ is a measurable positive function and 
     $$
     \lim_{x \to 0^+} \frac{h(tx)}{h(x)}=t^p \qquad  \forall t>0
     $$
     then one says that $h$ is regularly varying at the origin, in
     symbols $h \in R_p(0^+)$.
 \end{Def}
 
 Let $\tilde{h}(x):=h(\frac{1}{x})$.  Then $h \in R_p(0^+)$ iff
 $\tilde{h} \in R_{-p}$.
 
 \begin{Cor}\label{coro}
     $h \in R_1(0^+)\Longrightarrow \lim_{x \to 0^+}h(x)=0$.
 \end{Cor}
 \begin{proof}
     $h \in R_1(0^+)\Longrightarrow \tilde{h} \in R_{-1}$.  Therefore
     there exists $\ell$ slowly varying s.t.
     $\tilde{h}(x)=x^{-1}\ell(x)$.  This implies
     $$
     \lim_{x \to 0^+}h(x)=\lim_{y \to +\infty} h \left( \frac{1}{y}
     \right) =\lim_{y \to +\infty}\tilde{h}(y)= \lim_{y \to +\infty}
     \frac{\ell(y)}{y}=0
     $$
     where the last equality depends on Proposition \ref{eq:limiti}.
 \end{proof}

\section{The main result} \label{sect.V}

 \begin{Def}
     Two elements of $\gP$, $(\f,\c), (\tilde {\f}, \tilde {\c})$ are
     equivalent if there exist constants $A_1,A_2,B_1,B_2$ such that
     $A_1A_2=1$
     \begin{align*}
	 \tilde {\f} & = A_1\f+B_1 \\
	 \tilde {\c} & = A_2\c+B_2.
     \end{align*}
 \end{Def}
 
 Obviously equivalent elements of $\gP$ define the same $CM$-function.  In what
 follows we consider elements of $\gP$ up to this equivalence relation with
 the traditional abuse of language.
 
 \begin{Lem}\label{eq:lemma}
     Suppose that $(\f,\c)$ induce a paired monotone metric. Then 
     $$
     \frac{\f'(x)}{\f(x)}=\frac{p}{x} \quad  \quad \Longrightarrow \quad
     (\f,\c)=(\f_{p},\c_{p}).
     $$
 \end{Lem}
 \begin{proof}
     $$
     \frac{\f'(x)}{\f(x)}=\frac{p}{x}\Longrightarrow \log
     \frac{\f(x)}{\f(x_0)}=p \log \frac{x}{x_0}\Longrightarrow
     \f(x)=\f(x_0)(\frac{x}{x_0})^{p}=Ax^{p}.
     $$
     We may choose $\f(x)=\frac{x^p}{p}$.

     Now let $c(\cdot,\cdot)$ be the associated $CM$-function.  Going to
     the limit $y \to x$ in equation (\ref{eq:dualpair}) and using Proposition
     \ref{cm} one has
     \begin{equation}\label{derivata}
	 \f'(x)\c'(x)=c(x)=\frac{1}{x} 
     \end{equation}
     If $p=1$ then
     $$
     \f(x)=x \Longrightarrow \c'(x)=\frac{1}{x} \Longrightarrow
     \c(x)=\hbox{log}(x)
     $$
     If $p\not=1$ then
     $$
     \f(x)=x^p \Longrightarrow \c'(x)=\frac{1}{x^p} \Longrightarrow
     \c(x)= \frac{x^{1-p}}{1-p}
     $$
     and this ends the proof.
 \end{proof}
 
 We are ready to prove the fundamental result of the theory.

 \begin{Thm} \label{main} {\rm (Hasegawa and Petz (1997), Hasegawa 
 (2003))}
     Let $\f,\c\in\cc^{1}(0,+\infty)$.  Then $(\f,\c)$ induce a paired 
     monotone metric
     if and only if one of the following two possibilities holds
     $$
     (\f(x),\c(x))=(\frac{x^p}{p},\frac{x^{1-p}}{1-p}) \qquad 
     p\in[-1,2]\setminus\{0,1\}
     $$
     $$
     (\f(x),\c(x))=(x,\hbox{\rm{log}}(x)).
     $$
 \end{Thm}
 \begin{proof}
 
     The ``if" part is just Theorem \ref{family}.  To prove the ``only
     if" part we need some auxiliary functions.

     Let us  define
     $$
     k(x,y):= (\f(x)-\f(y))(\c(x)-\c(y))=(x-y)^2c(x,y).
     $$
     One has
     $$
     k(tx,ty)=t^2(x-y)^2c(tx,ty)=t(x-y)^2c(x,y)=tk(x,y)
     $$
     that is $k$ is 1-homogeneous. Moreover set $h(x):=\f(x)\c(x)$.

     Equation (\ref{derivata}) implies that $\f,\c$ are strictly
     monotone (either both increasing or both decreasing) and
     therefore injective.  Moreover monotonicity implies that the
     following limits exist
     $$
     \f(0^+):=\lim_{x \to 0+}\f(x) \qquad \qquad \c(0^+):=\lim_{x \to
     0+}\c(x).
     $$
     Since we consider $\f,\c$ up to additive constants and because
     we can change the sign of $\f,\c$, we may reduce to $\f,\c$ 
     increasing, and  have to consider three cases

     a) $\f(0^+)=\c(0^+)=-\infty $,
     
     b) $\f(0^+)=0,\  \c(0^+)=-\infty $,
     
     c) $\f(0^+)=\c(0^+)=0$.

     \bigskip
     
     \noindent Case a)
     \bigskip
     
     Suppose $\f(0^+)=\c(0^+)=-\infty$.  Now let $0<y<x$; going to the
     limit $y \to 0^+$ we have that
     \begin{align*}
	 t & = \frac{k(tx,ty)}{k(x,y)}= \lim_{y \to 0^+}
	 \frac{\f(tx)-\f(ty)}{\f(x)-\f(y)}\cdot\frac{\c(tx)-\c(ty)}{\c(x)-\c(y)}
	 \\
	 & =\lim_{y \to
	 0^+}\frac{\f(ty)}{\f(y)}\cdot\frac{\c(ty)}{\c(y)}=\lim_{y \to
	 0^+}\frac{h(ty)}{h(y)}.
     \end{align*}
     This means that $h \in R_1(0^+)$ and therefore by Corollary
     \ref{coro}
     $$
     +\infty=\lim_{x \to 0^+}\f(x)\c(x)=\lim_{x \to 0^+}h(x)=0
     $$
     that is absurd.
     \bigskip
     
     \noindent Case b)
     \bigskip

     Suppose $\f(0^+)=0$ and $\c(0^+)=-\infty$.  Again let $0<y<x$;
     going to the limit $y \to 0^+$ we have that
     $$
     t = \lim_{y \to 0^+} \frac{\f(tx)-\f(ty)}{\f(x)-\f(y)} \cdot
     \frac{\c(tx)-\c(ty)}{\c(x)-\c(y)} = \frac{\f(tx)}{\f(x)} \cdot
     \lim_{y \to 0^+} \frac{\c(ty)}{\c(y)}.
     $$
     This implies that the limit
     $$
     \lim_{y \to 0^+}\frac{\c(ty)}{\c(y)}
     $$
     exists $\forall t$.  Therefore there exists a function $j$ such
     that
     $$
     j(t)=\lim_{y \to 0^+}\frac{\c(ty)}{\c(y)}=\frac{t\f(x)}{\f(tx)}.
     $$
     From Proposition \ref{RV} one has that $-\c \in R_\b(0^+)$ for
     some $\b \in \br$ namely $j(t)=t^\b$, $\forall t>0$.  So we have
     $\f(tx)=t^{1-\b}\f(x)$, that is $\f$ is $p$-homogeneous, with
     $p:= 1-\b$, and therefore by Euler $x\f'(x)=p\f(x)$ ($p \not= 0$
     because $\f'\not=0$); then
     $$
     \frac{\f'(x)}{\f(x)}=\frac{p}{x}
     $$
     and therefore because of Lemma \ref{eq:lemma},
     $$
     (\f,\c) = (\f_{p},\c_{p}).
     $$
     Since $\c(0^+)=-\infty$ we have $p\geq1$.  Because of Theorem
     \ref{family} one has $p\in[1,2]$.
     \bigskip
     
     \noindent Case c)
     \bigskip

     The argument for this case is that of Hasegawa and Petz (1997),
     Hasegawa (2003) and we report it here for the sake of
     completeness.

     One can deduce
     $$
     \frac{\f(tx)-\f(ty)}{t(x-y)}\cdot
     \frac{\c(tx)-\c(ty)}{t(x-y)}=\frac{1}{t}\frac{\f(x)-\f(y)}{x-y}\cdot
     \frac{\c(x)-\c(y)}{x-y}
     $$
     Going to the limit $y \to 0^+$ one has
     $$
     \frac{\f(tx)}{tx}\cdot
     \frac{\c(tx)}{tx}=\frac{1}{t}\frac{\f(x)}{x}\cdot \frac{\c(x)}{x}
     $$ 
     so that
     $$
     \f(tx)\c(tx)=t\f(x)\c(x)
     $$
     This means that $h(x)=\f(x)\c(x)$ is 1-homogeneous and $h(0^+)=0$
     so that, because of Euler, one has $xh'(x)=h(x)$.  This implies
     that $\exists b \in \br$ s.t. $h(x)=bx$, $\forall x>0$.  Then
     $$
     \f(x)\c(x)=bx \qquad \qquad   \f'(x)\c'(x)= \frac{1}{x}.
     $$
     As $\f,\c$ are increasing, $b>0$.  Deriving the first equation
     one gets
     $$
     \f'(x)\c(x)+\f(x)\c'(x)=b.
     $$
     Since $\f,\f'\not= 0$, $\forall x>0$ one may write 
     $$
     \c(x)=\frac{bx}{\f(x)} \qquad \qquad \c'(x)=\frac{1}{x\f'(x)}.
     $$ 
     Substituting one gets
     $$
     \frac{\f'(x)}{\f(x)} bx+ \frac{\f(x)}{\f'(x)} \frac{1}{x}=b,
     $$
     so that if $y(x):=\frac{\f'(x)}{\f(x)}\not= 0$ the equation
     becomes
     $$
     bxy(x)+ \frac{1}{xy(x)}=b
     $$
     and finally
     $$
     bx^2y(x)^2-bxy(x)+1=0.
     $$
     From this it follows that

     i) if $0 < b <4$  there is no solution;

     ii) if $b \geq 4$ then 
     $$
     y(x)=\frac{1}{2x}\left(1\pm\sqrt{1-\frac{4}{b}}\right).
     $$
     Therefore, setting 
     $$p:=\frac{1+ \sqrt{1-\frac{4}{b}}}{2}\in\left[\frac12,1\right),$$
     we have
     $$\frac{\f'(x)}{\f(x)}= y(x)= \frac{p}{x}$$
     or
     $$\frac{\f'(x)}{\f(x)}= y(x)= \frac{1-p}{x}.$$
     From Lemma \ref{eq:lemma} one has
     $$
     \f(x)=\frac{x^p}{p} \qquad \c(x)=\frac{x^{1-p}}{1-p}
     $$
     or viceversa.  Therefore $(\f,\c) = (\f_{p},\c_{p})$, with
     $0<p<1$, and this ends the proof.
 \end{proof}

 \begin{Cor}
     If $\f,\c$ induce a paired monotone metric then $\lim_{x \to
     0+}\f(x)\c(x)=0$.
 \end{Cor}

 \begin{Rem} \label{condizioni}
     In Hasegawa and Petz (1997), Hasegawa (2003) Theorem \ref{main}
     is proved under the hypothesis that $\lim_{x \to
     0+}\f(x)\c(x)=0$.  The present proof shows that this hypothesis
     can be dropped.  An application of this is given in Gibilisco and 
     Isola (2003).
 \end{Rem}

 \section{Concluding remarks} 
 
 In Hasegawa (2003) Hasegawa wanted to find a family of (non-paired) operator
 monotone functions that "fill the gap" between the functions
 $$
 f_{\rm Bures}(x)=\frac{1+x}{2} \qquad \qquad
 f_{\frac{1}{2}}(x)=\frac{(1+ \sqrt{x})^2}{4}
 $$ 
 corresponding to the SLD-metric and the WY-metric.  The problem can
 be solved by proving the following 
 
 \begin{Prop}\label{gap}
     The functions
     $$
     f^{\nu}_{\rm power}(x):=\left( \frac{1+x^{\frac{1}{\nu}} }{2}
     \right)^{\nu} \qquad \qquad 1 \leq \nu \leq 2
     $$
     are operator monotone.
 \end{Prop}

 To prove Proposition \ref{gap} Hasegawa used an argument due to Petz. 
 We just want to remark that the above result can be proved by applying
 to $f_{\rm Bures}(x)=\frac{1+x}{2}$ the following
 
 \begin{Prop}\label{ando}
     Let $f$ be operator monotone, and $\n\in[1,\infty)$.  Then
     $x\in(0,\infty) \to f(x^{1/\n})^\n$ is operator monotone.
 \end{Prop}
 \begin{proof}
      See the proof of Corollary 4.3 $(i)$ in Ando (1979).
 \end{proof}

 For $p\in(0,1)$ namely $q\in(1,+\infty)$ Proposition \ref{Banach}
 shows that it is possible to relate the duality discussed here to the
 geometry of spheres in $L^q$ spaces along the lines of Gibilisco and
 Pistone (1998), Gibilisco and Isola (1999, 2001a,b).  The same does
 not apply to the cases $p\in[-1,0]$ or $p=1$.  In Gibilisco and
 Pistone (1998) the Amari embedding was generalised to the sphere of
 an Orlicz space under very general hypothesis.  We conjecture that,
 for $p=0,1$ (that is for the BKM metric) one can use non-commutative
 analogues of the Zygmund spaces $L^{{\rm exp}}$, $L^{x{\rm log}x}$ to
 produce a similar construction (see also Grasselli and Streater
 (2000) and references therein).

 \section*{Acknowledgements} 
 It is a pleasure to thank H.Hasegawa for useful conversations on this
 subject at the congress "Information Geometry and its Applications",
 Pescara, Italy, 2002.  This research has been supported by the
 italian MIUR program "Quantum Probability and Infinite Dimensional
 Analysis" 2001-2002.

%%%%%%REFERENCES%%%%%%%%
 \section*{References}

 \quad Amari, S. and Nagaoka, H. (2000). {\it Methods of Information
 Geometry},  American Mathematical Society and Oxford University
 Press.

 Ando, T. (1979).  Concavity of certain maps on positive definite
 matrices and applications to Hadamard products, {\it Linear Algebra
 and its Applications}, {\bf 26}, 203--241.
 
 Berger, M. S. (1977). {\it Nonlinearity and functional 
 analysis}. Academic Press, New York.
 
 Bhatia, R. (1997). {\it Matrix Analysis}. Springer-Verlag,
 New York.

 Bingham, N.H., Goldie, C.M. and Teugels, J.L. (1987).
 {\it Regular variation}, Cambridge University Press, Cambridge.

 Chentsov, N. (1982). {\it Statistical decision rules and
 optimal inference}, American Mathematical Society, Providence.

 Chentsov, N. and Morotzova, E. (1990). Markov invariant geometry
 on state manifolds (in russian), {\it Itogi Nauki i Tekhniki}, {\bf
 36}, 69-102.
 
 Gibilisco, P. and Isola, T. (1999). Connections on statistical 
 manifolds of density operators by geometry of noncommutative 
 $L^{p}$-spaces,  {\it Infinite Dimensional Analysis, Quantum Probability \& 
 Related Topics},  {\bf 2}, 169-178.
 \smallskip
 
 Gibilisco, P. and Isola, T. (2001a). Monotone metrics on
 statistical manifolds of density matrices by geometry of
 noncommutative $L^{2}$-spaces, {\it Disordered and complex
 systems}, 129-140, Amer.  Inst.  Phys., Melville.  

 Gibilisco, P. and Isola, T. (2001b).  A characterisation
 of Wigner-Yanase skew information among statistically monotone
 metrics,  {\it Infinite Dimensional Analysis, Quantum Probability \& 
 Related Topics},   {\bf 4}, 553-557. 

 Gibilisco, P. and Isola, T. (2003).  Wigner-Yanase information on
 quantum state space: the geometric approach, {\it Journal of 
 Mathematical  Physics}, 
 {\bf 44}, 3752-3762.

 Gibilisco, P. and Pistone G. (1998).  Connections on non-parametric
 statistical manifolds by Orlicz space geometry, {\it Infinite
 Dimensional Analysis, Quantum Probability \& Related Topics}, {\bf
 1}, 325-347.

 Grasselli, M.R. (2002). Duality, monotonicity and the 
 Wigner-Yanase-Dyson metrics, preprint 2002, [math-ph/0212022].

 Grasselli, M.R. and Streater R.F. (2000).  The quantum information
 manifold for $\epsilon$-boundend forms, {\it Reports on Mathematical
 Physics }, {\bf 46}, 325-335.

 Grasselli, M.R. and Streater R.F. (2001).  On the uniqueness of
 Chentsov metric in quantum information geometry, {\it Infinite
 Dimensional Analysis, Quantum Probability \& Related Topics}, {\bf
 4}, 173-182.
 
 Hasegawa, H. (1993).  $\alpha$-divergence of the non-commutative
 information geometry, {\it Reports on Mathematical Physics}, {\bf 33},
 87-93.

 Hasegawa, H. (1995). Non-commutative extension of the
 information geometry, {\it Quantum Communications and Measurement}, 
 327-337, Plenum Press, New York.

 Hasegawa, H. and Petz, D. (1997).  Non-commutative extension of
 information geometry, II, {\it Quantum Communication, Computing and
 Measurement}, 109-118, Plenum Press, New York.

 Hasegawa, H. (2003).  Dual Geometry of the Wigner-Yanase-Dyson
 Information Content.  Preprint.  To appear on {\it Infinite
 Dimensional Analysis, Quantum Probability \& Related Topics}.

 Jen\v{c}ova, A. (2001).  Geometry of quantum states: dual connections
 and divergence functions, {\it Reports on Mathematical Physics}, {\bf
 47}, 121-138.
 
 Lesniewski, A. and Ruskai, M.~B. (1999).  Monotone Riemannian metrics
 and relative entropy on noncommutative probability spaces, {\it
 Journal of Mathematical Physics}, {\bf 40}, 5702--5724.
 
 Lieb E.H. (1973). Convex trace functions and the
 Wigner-Yanase-Dyson conjecture, {\it Advances in Mathematics}, {\bf 11},
 267-288.

 Nagaoka H. (1995).  Differential geometric aspects of quantum state
 estimation and relative entropy, {\it Quantum Communications and
 Measurement}, 449-452, Plenum Press, New York.

 Petz, D. (1996). Monotone metrics on matrix space, 
 {\it Linear Algebra and its Applications},  {\bf 244}, 81--96. 
 
 Petz, D. (2002).  Covariance and Fisher information in
 quantum mechanics, {\it Journal of Physics. A. Mathematical and General}, {\bf 35}, 929-939.  
 
 Petz, D. and Hasegawa, H. (1996).  On the Riemannian metric of 
 $\alpha$-entropies of density matrices, {\it Letters in Mathematical Physics},  
 {\bf 38}, 221-225.
 
 Pistone G. and Sempi C. (1995). An infinite-dimensional geometric
 structure on the space of all probability measures equivalent to a
 given one, {\it Annals of  Statistics}, {\bf 23}, 1543-1561.
 
 Wigner, E. and Yanase M. (1963).  Information content of
 distribution, {\it Proceedings of the National Academy of Sciences of
 the United States of America}, {\bf 49}, 910-918.

\end{document}